\documentclass{article}

\usepackage[final]{nips_2017}

\usepackage[utf8]{inputenc} 
\usepackage[T1]{fontenc}    
\usepackage{hyperref}       
\usepackage{url}            
\usepackage{booktabs}       
\usepackage{amsfonts}       
\usepackage{nicefrac}       
\usepackage{microtype}      

\usepackage{amsmath}
\usepackage{dsfont}
\usepackage{varioref}

\title{Dimension-free PAC-Bayesian bounds for the estimation of the mean 
of a random vector}

\author{
Olivier Catoni 
\\
CREST -- CNRS UMR 9194 \\
Université Paris Saclay \\ 
  \texttt{olivier.catoni@ensae.fr} 
  \And
Ilaria Giulini\\
Laboratoire de Probabilités\\ et Modèles Aléatoires\\
Université Paris Diderot\\
\texttt{giulini@math.univ-paris-diderot.fr} 
}

\newtheorem{thm}{Theorem}[section]
\newtheorem{prop}[thm]{Proposition}

\newtheorem{lemma}[thm]{Lemma}

\newcommand{\ds}{\displaystyle}

\renewcommand{\mathbb}{\mathds}

\newcommand{\C}[1]{\mathcal{#1}}

\newcommand{\wt}[1]{\widetilde{#1}}
\newcommand{\wh}[1]{\widehat{#1}}
\newcommand{\B}[1]{\mathds{#1}}

\newcommand{\ud}{\mathrm{d}}

\begin{document}

\maketitle

\begin{abstract}
	In this paper, we present a new estimator of the mean of a 
random vector, computed by applying some threshold function 
to the norm. Non asymptotic dimension-free almost sub-Gaussian 
bounds are proved under weak moment assumptions, using PAC-Bayesian 
inequalities.  
\end{abstract}

\section{Introduction}

Estimating the mean of a random vector under weak tail assumptions
has attracted a lot of attention recently. A number of properties 
have spurred the interest for these new results, where the empirical 
mean is replaced by a more robust estimator. One aspect is 
that it is possible to obtain an estimator with a sub-Gaussian 
tail while assuming much weaker assumptions on the data, 
up to the fact of assuming only the existence 
of a finite covariance matrix. Another appealing feature is 
that it is possible to obtain dimension-free non asymptotic 
bounds that remain valid in a separable Hilbert space.
Some important references are \citet{Cat10} in the one dimensional 
case and \citet{Minsker} and \citet{LugoMen} in the multidimensional 
case. Building on the breakthrough of \citet{Minsker}, 
that uses a multidimensional generalization of the median of means 
estimator, \citet{JolyLugoOl} and \citet{LugoMen} propose 
successive improvements of the median of means approach
to get an estimator with a genuine sub-Gaussian 
dimension-free tail bound, 
while still requiring only the existence of the covariance matrix.
In the mean time, the M-estimator approach of \citet{Cat10} 
has also been generalized to multidimensional settings through 
the use of matrix inequalities in \citet{Minsker2} and \citet{MinskerWei}.

Here we follow a different route, based on a multidimensional 
extension of \citet{Cat10} using PAC-Bayesian bounds. Our new estimator 
is a simple modification of the empirical mean, where some 
threshold is applied to the norm of the sample vectors. 
Therefore, it is straightforward to compute, and this is
a strong point of our approach, compared to others.
Note also that we make here some compromise on the 
sharpness of the estimation error bound, in order to 
simplify the definition and computation of the 
estimator. This compromise consists in the presence of 
second order terms, while the first order terms can 
be made as close as desired to a true sub-Gaussian bound 
with exact constants, as stated in \citet[eq. (1.1)]{LugoMen}. 
With a more involved estimator, a true 
sub-Gaussian bound without second order terms is 
possible and will be described in a separate 
publication. 

\section{Thresholding the norm} 

Consider $X \in \B{R}^d$, a random vector, 
and $(X_1, \dots, X_n)$ a sample made of $n$ independent 
copies of $X$. The question is to estimate 
$\B{E}(X)$ from the sample, under the assumption that 
$\B{E} \bigl( \lVert X \rVert^p \bigr) < \infty$, for some $p \geq 2$. 
 
Consider the threshold function $\ds 
\psi(t) = \min \{ t, 1 \}, \; t \in \B{R}_+$,
and for some positive real parameter $\lambda$ 
to be chosen later, introduce the thresholded 
sample 
\[
Y_i = \frac{\psi \bigl( \lambda \lVert X_i \rVert \bigr)}{ \lambda 
\lVert X_i \rVert} X_i. 
\]
Our estimator of $m = \B{E}(X)$ will simply be the thresholded empirical mean 
$\ds \wh{m} = \frac{1}{n} \sum_{i=1}^n Y_i$.
\begin{prop}
\label{prop:2.1}
Introduce the increasing functions 
\[
g_1(t) = \frac{1}{t} \Bigl(\exp(t) - 1 \Bigr)
\text{ and } 
g_2(t) = \frac{2}{t^2} \Bigl( \exp(t) - 1 - t \Bigr), \qquad t \in \B{R},
\]
that are defined by continuity at $t = 0$ and are such that $g_1(0) = g_2(0) = 
1$. 
Assume that $\B{E} \bigl( \lVert X \rVert^2 \bigr) < \infty$ and that 
we know $v$ such that  
\[
\sup_{\theta \in \B{S}_d} \B{E} \bigl( \langle \theta, X - m \rangle^2 
\bigr) \leq v  < \infty, \\ 
\]
where $\ds \B{S}_d = \bigl\{ \theta \in \B{R}^d, \lVert \theta \rVert = 1 
\bigr\}$ is the unit sphere of $\B{R}^d$. 
For some positive real parameter $\mu$, put 
\begin{align*}
\lambda & = \mu^{-1} \sqrt{\frac{2 \log ( \delta^{-1})}{a v n}}, &  
T & = \max \bigl\{ \B{E} \bigl( \lVert X - m \rVert^2 \bigr), v \bigr\}, \\ 
a & = g_2(2 \mu) \geq 1, & 
b & = \exp(2 \mu) g_1 \Biggl( \mu^2 \sqrt{\frac{2 a v}{ 
T \log(\delta^{-1})}} \; \Biggr) \geq 1.  
\end{align*}
With probability at least $ 1 - \delta$,
\[
\lVert \wh{m} - m \rVert \leq \sqrt{ \frac{2 a v \log(\delta^{-1})}{n}} + 
\sqrt{\frac{b T}{n}} + \inf_{p \geq 1} \frac{C_p}{ n^{p/2}} + \inf_{p \geq 2} 
\frac{C'_p}{n^{p/2}},
\]
where 
\begin{align*}
C_p & = \frac{1}{p+1} \biggl( 
\frac{p}{(p+1) \mu } \biggr)^p \biggl( \frac{ 2 \log(\delta^{-1})}{a v} 
\biggr)^{p/2}  \sup_{\theta \in \B{S}_d} 
\B{E} \bigl( \lVert X \rVert^p \langle \theta, X - m \rangle_- 
\bigr), \quad \text{ and } \\ 
C'_p & = \frac{1}{p+1} \biggl( \frac{p}{(p+1) \mu} \biggr)^p \biggl( \frac{
2 \log(\delta^{-1})}{a v} \biggr)^{p/2} \B{E} \bigl( 
\lVert X \rVert^p \bigr) \lVert m \rVert 
\Biggl( 1 + 
\sqrt{\frac{a \log(\delta^{-1})}{2 v n}}
\lVert m \rVert \Biggr). 
\end{align*}
\end{prop}
\paragraph{Remarks}
Note that in case $\B{E} \bigl( \lVert X \rVert^2 \bigr) < \infty$ 
but $\B{E} \bigl( \lVert X \rVert^p \bigr)  = \infty$ for $p > 2$, 
we can use the bound 
\begin{multline*}
\frac{C_1}{\sqrt{n}} + \frac{C'_2}{n} \leq 
\frac{1}{2 \mu} \sqrt{\frac{\ds  \log(\delta^{-1})( T + \lVert m \rVert^2)}{ 
2 a n }} 
+ \frac{8 \log(\delta^{-1})}{27 \mu^2 a v n} \B{E} \bigl( 
\lVert X \rVert^2 \bigr)  \lVert m \rVert \\ \times \Biggl( 
1 + \sqrt{\frac{ a \log(\delta^{-1})}{2 v n}} 
 \lVert m \rVert \Biggr) = 
\C{O} \Biggl( \frac{1}{2 \mu} \sqrt{\frac{\log(\delta^{-1}) ( T + 
\lVert m \rVert^2)}{2 a n}} \; \Biggr).
\end{multline*}
Note also that if we take $\mu = 1/4$ and assume that $\delta \leq \exp( - 1)$, 
then $a \leq 1.2$ and $b \leq 4$. 
If moreover 
$\B{E} \bigl( \lVert X \rVert^{p+1} 
\bigr) < \infty$, for some $p > 1$, 
we obtain with probability at least $1 - \delta$ that 
\[ 
\lVert \wh{m} - m \rVert \leq \sqrt{ \frac{2.4 
\, v \log( 
\delta^{-1})}{n}} + \sqrt{\frac{4 T}{n}}
+ \frac{C_p}{n^{p/2}} + \frac{C'_{p+1}}{n^{(p+1)/2}},  
\] 
meaning that the tail distribution of $\lVert \wh{m} - m \rVert$
has a sub-Gaussian behavior, up to second order terms. 
Remark that by taking $\mu$ small, we can make $a$ and $b$ 
as close as desired to $1$, at the expense of the values 
of $C_p$ and $C'_p$. 

\paragraph{Proof} The rest of the paper is devoted to the proof of 
Proposition \ref{prop:2.1}. 

An elementary computation shows that the threshold function 
$\psi$ satisfies 
\begin{equation}
\label{eq:1.2}
0 \leq 1 - \frac{\psi(t)}{t} \leq \inf_{p \geq 1} \frac{t^{p}}{p + 1}
\biggl( \frac{p}{p+1} \biggr)^p, \qquad t \in \B{R}_+, 
\end{equation}
where non integer values of the exponent $p$ are allowed. 
Let $\ds Y = \frac{\psi \bigl( \lambda \lVert X \rVert \bigr)}{\lambda 
\lVert X \rVert} X$ and $\wt{m} = \B{E}(Y)$. We can decompose 
the estimation error in direction $\theta$ into 
\begin{equation}
\label{eq:1}
\langle \theta, \wh{m} - m \rangle = \langle \theta, \wt{m} - m \rangle 
+ \frac{1}{n} \sum_{i=1}^n \langle \theta, Y_i - \wt{m} \rangle, \qquad 
\theta \in \B{R}^d.
\end{equation}
Introduce $\ds \alpha = \frac{\psi \bigl( 
\lambda \lVert X \rVert \bigr) }{ \lambda \lVert X \rVert}$ 
and let us deal with the first term first. As $\ds 0 \leq 1 - \alpha 
\leq \frac{\lambda^p \lVert X \rVert^p}{p+1} \biggl( \frac{p}{p+1} \biggr)^p$
\begin{multline*}
\langle \theta, \wt{m} - m \rangle = \B{E} 
\bigl[ ( \alpha - 1 ) \langle \theta, X \rangle \bigr] 
= \B{E} \bigl[ ( \alpha - 1 ) \langle \theta, X - m \rangle
\bigr]  + \B{E}  
( \alpha - 1 ) \langle \theta, m \rangle 
\\ \leq \inf_{p \geq 1} 
\frac{\lambda^p}{(p+1)} \biggl( \frac{p}{p+1} 
\biggr)^p 
\B{E} \Bigl( \lVert X \rVert^p \langle \theta, X - m  \rangle_- 
\Bigr) 
+ \inf_{p \geq 2} \frac{\lambda^p}{(p+1)} \biggl( \frac{p}{p+1} 
\biggr)^p \B{E} \bigl( \lVert X \rVert^p \bigr) \langle \theta, 
m \rangle_-,
\end{multline*}
where $r_- = \max \{ 0, -r \}$ is the negative part of integer 
$r$. 

Let us now look at the second term of the decomposition \eqref{eq:1}. 
To gain uniformity in $\theta$, we will use a PAC-Bayesian 
inequality and the family of normal distributions 
$\rho_{\theta} = \C{N} \bigl(\theta, 
\beta^{-1} I_d \bigr)$, bearing on the parameter $\theta \in \B{R}^d$, 
where $I_d \in \B{R}^{d \times d}$ is the identity matrix of size $d 
\times d$, and where $\beta$ is a positive parameter to be chosen 
later on.   

We will use the following PAC-Bayesian inequality without recalling 
its proof, that is a simple consequence of \citet[eq. (5.2.1) page 159]{Cat01b}:
\begin{lemma}
\label{lem:2.2}
For any bounded measurable function $f : \B{R}^d \times \B{R}^d \rightarrow 
\B{R}$, for any probability measure $\pi \in \C{M}_+^1 \bigl(
\B{R}^d \bigr)$, for any $\delta \in ]0,1[$, with probability at least $1 - \delta$, 
for any probability measure $\rho \in \C{M}_+^1(\B{R}^d)$, 
\[
\frac{1}{n} \sum_{i=1}^n \int f \bigl( \theta, X_i \bigr) \, \ud 
\rho(\theta) \leq \int \log \Bigl[ \B{E} \Bigl( \exp \bigl( 
f(\theta, X) \bigr) \Bigr) \Bigr] \, \ud \rho(\theta) \\ 
+ \frac{\C{K}(\rho, \pi) + \log(\delta^{-1})}{n}, 
\]
where $\C{K}$ is the Kullback-Liebler divergence 
$\ds
\C{K} (\rho, \pi) = \begin{cases}
\int \log \bigl( \rho / \pi \bigr) \, \ud \rho, & \text{ when } \rho \ll \pi, \\
+ \infty, & \text{ otherwise.}
\end{cases}
$
\end{lemma} 
Remarking that 
$\ds
\frac{1}{n} \sum_{i=1}^n \langle \theta, Y_i - \wt{m} \rangle = 
\frac{1}{n} \sum_{i=1}^n \int \langle \theta', Y_i - \wt{m} \rangle 
\, \ud \rho_{\theta}( \theta' ),  
$
using $\pi = \rho_0$, and taking into account the fact that  
$\C{K}(\rho_{\theta}, \rho_0 ) = \beta \lVert \theta \rVert^2 / 2$, 
we obtain as a consequence of the previous lemma that 
with probability at least $1 - \delta$, for any 
$\theta \in \B{S}_d$, 
\[
\frac{1}{n} \sum_{i=1}^n \langle \theta, Y_i - \wt{m} \rangle 
\leq \frac{1}{\mu \lambda} \int \log \biggl( \B{E} 
\exp \Bigl( \mu \lambda \langle \theta', Y - \wt{m} \rangle \Bigr) \biggr)
\, \ud \rho_{\theta}(\theta')
+ \frac{\beta}{2 n \mu \lambda} + \frac{\log ( \delta^{-1})}{ n \mu \lambda}. 
\]
In our setting $f$ is not bounded in $\theta$, but the required extension 
is valid as explained in \citet{Cat01b}. Since the logarithm is concave, 
\begin{multline*}
\int \log \biggl( \B{E} \exp \Bigl( \mu \lambda \langle \theta', 
Y - \wt{m} \rangle \Bigr) \biggr) \, \ud \rho_{\theta}(\theta') 
\leq \log \Biggl[ \B{E} \biggl( 
\int \exp \Bigl( \mu \lambda \langle \theta', 
Y - \wt{m} \rangle \Bigr) \, \ud \rho_{\theta}(\theta') \biggr) \Biggr]
\\ = 
\log \Biggl[ \B{E} \biggl( \exp \Bigl( \mu \lambda \langle \theta, 
Y - \wt{m} \rangle + \frac{\mu^2 \lambda^2}{2 \beta} \lVert Y - \wt{m} \rVert^2 
\Bigr) \biggr) \Biggr],
\end{multline*}
where we have used the explicit expression of the Laplace transform 
of a Gaussian distribution.

To go further, reminding as a source of inspiration the proof of 
Bennett's inequality, let us introduce the increasing functions $g_1$ and $g_2$
defined in Proposition \ref{prop:2.1}.
These functions will be used to bound the exponential 
function by polynomials. More precisely, we will exploit the fact 
that when $t \leq b$, $\exp(t) \leq 1 + t + g_2(b) t^2 / 2$ and 
$\exp(t) \leq 1 + g_1(b) t$. From this, it results that if 
$t \leq b$ and $u \leq c$, 
\begin{multline*}
\exp( t + u) \leq \exp(t) \bigl(1 + g_1(c) u \bigr) 
\leq \exp(t) + g_1(c) \exp(b) u \\ \leq 
1 + t + g_2(b) t^2/2 + g_1(c) \exp(b) u.
\end{multline*}

Legitimate values for $b$ and $c$ will be deduced 
from the remark that $\lambda \lVert Y \rVert \leq 1$, implying 
$\lambda \lVert \wt{m} \rVert \leq 1$. Namely, in our context, 
we will use $b = 2 \mu$ and $c = 2 \mu^2 / \beta$. 

These arguments put together lead to the inequality
\begin{multline*}
\B{E} \Biggl( \exp \biggl( \mu \lambda 
\langle \theta, Y - \wt{m} \rangle + \frac{ \mu^2 \lambda^2}{2 
\beta} \bigl\lVert Y - \wt{m} \bigr\rVert^2 \biggr) \Biggr) 
\\ \leq 1 + g_2(2 \mu) \frac{\mu^2 \lambda^2}{2} \B{E} \bigl( \langle \theta, 
Y - \wt{m} \rangle^2 \bigr) + \exp (2 \mu ) g_1 \biggl( 
\frac{2 \mu^2}{\beta} \biggr) \frac{\mu^2 \lambda^2}{2 \beta} 
\B{E} \Bigl( \bigl\lVert Y - \wt{m} \bigr\rVert^2 \Bigr).
\end{multline*}

Replacing in the previous inequalities, 
we obtain
\begin{lemma}
With probability at least $1 - \delta$, for any $\theta \in \B{S}_d$, 
\begin{multline*}
\langle \theta, \wh{m} - \wt{m} \rangle = \frac{1}{n} \sum_{i=1}^n \langle \theta, Y_i - \wt{m} 
\rangle \leq g_2(2 \mu) \frac{\mu \lambda}{2} \B{E} \bigl( 
\langle \theta, Y - \wt{m} \rangle^2 \bigr) \\ 
+ \exp(2 \mu) g_1 \biggl( \frac{2 \mu^2}{\beta} \biggr) 
\frac{\mu \lambda}{2 \beta} 
\B{E} \bigl( \lVert Y - \wt{m} \rVert^2 \bigr) + 
\frac{ \beta + 2 \log(\delta^{-1})}{2 \mu \lambda n}.
\end{multline*}
\end{lemma}
Remark that
\begin{multline*}
\langle \theta, Y - m \rangle^2  
= \langle \theta, \alpha X - m \rangle^2 
= \Bigl( \alpha \langle \theta, X - m \rangle - (1 - \alpha) \langle 
\theta, m \rangle \Bigr)^2
\\ \leq \alpha \langle \theta, X - m \rangle^2 + (1 - \alpha) \langle \theta, m 
\rangle^2 \leq \langle \theta, X - m \rangle^2 + (1 - \alpha) \langle \theta, 
m \rangle^2. 
\end{multline*}
Therefore, using inequality \eqref{eq:1.2} and the definition of $\alpha$, 
\[
\B{E} \bigl( \langle \theta, Y - \wt{m} \rangle^2 \bigr) \leq 
\B{E} \bigl( \langle \theta, Y - m \rangle^2 \bigr) \leq 
\B{E} \bigl( \langle \theta, X - m \rangle^2 \bigr) + \langle \theta, 
m \rangle^2 \inf_{p \geq 2} \frac{\lambda^p}{p+1} \biggl( \frac{p}{p+1} 
\biggr)^p \B{E} \bigl( \lVert X \rVert^p \bigr).  
\]
Remark also that $Y = g(X)$, where $g$ is a contraction (being the 
projection on a ball). Consequently
\[
\B{E} \bigl( \lVert Y - \wt{m} \rVert^2 \bigr) = \frac{1}{2} 
\B{E} \bigl( \lVert Y_1 - Y_2 \rVert^2 \bigr) \leq  
\frac{1}{2} \B{E} \bigl( \lVert X_1 - X_2 \rVert^2 \bigr) = \B{E} 
\bigl( \lVert X - m \rVert^2 \bigr).  
\]
In view of these remarks, the previous lemma translates to
\begin{lemma}
Let $\ds a = g_2 \bigl( 2 \mu \bigr)$ and 
$\ds b \geq \exp ( 2 \mu ) g_1 \biggl( \frac{2 \mu^2}{\beta} \biggr)$.\\ 
With probability at least $1 - \delta$, for any $\theta \in \B{S}_d$, 
\begin{multline*}
\langle \theta, \wh{m} - m \rangle \leq 
\frac{a \mu \lambda}{2}  
\B{E} \bigl( \langle \theta, X - m \rangle^2 \bigr) 
+ \frac{b \mu \lambda}{2 \beta} \B{E} \bigl( \lVert X - m \rVert^2 \bigr) 
+ \frac{\beta + 2 \log(\delta^{-1})}{2 \mu \lambda n} 
\\ + 
\inf_{p \geq 1} 
\frac{\lambda^p}{p+1} \biggl( \frac{p}{p+1} \biggr)^p 
\B{E} \bigl( \lVert X \rVert^p \langle \theta, X - m\rangle_- \bigr)
\\ + \inf_{p \geq 2} \frac{\lambda^p}{p+1} \biggl( \frac{p}{p+1} 
\biggr)^p \B{E} \bigl( \lVert X \rVert^p \bigr) 
\Biggl( \langle \theta, m \rangle_- + \frac{a \mu \lambda}{2} \langle 
\theta, m \rangle^2 
\Biggr). 
\end{multline*}
\end{lemma}
Proposition \ref{prop:2.1} follows by taking $b$ as mentioned there, 
$\ds\lambda = \frac{1}{\mu} \sqrt{\frac{2 \log(\delta^{-1})}{n a v}}$, 
and $\ds \beta =  \sqrt{\frac{2 b T \log(\delta^{-1})}{a v}} \geq \sqrt{
\frac{2 T \log(\delta^{-1})}{av}}$, so that the condition on $b$ 
is satisfied. 

\small

\bibliography{CatGiu}

\end{document}